\documentclass[twocolumn]{autart}    
\usepackage{xcolor}\usepackage{amsmath}\usepackage{indentfirst}
        \usepackage{graphicx}
\usepackage{subfigure} 

\usepackage{amsmath}
\usepackage{amsfonts}
\usepackage{graphicx} 
\usepackage{subfigure}

\begin{document}

\begin{frontmatter}

\title{Maximum principle for discrete-time stochastic recursive games driven by fractional noises$^{\ast,\ast\ast\ast}$} 

\thanks[footnoteinfo]{This paper was not presented at any IFAC 
meeting. }
\thanks[footnoteinfo]{\quad Corresponding author: Chenfei Wang.}
\thanks[footnoteinfo]{\qquad 
  These authors contributed equally to this work.}
\author[Paestum]{Yuecai Han}\ead{hanyc@jlu.edu.cn},    
\author[Paestum]{Yuhang Li}\ead{yuhangl22@mails.jlu.edu.cn}, 
\author[Paestum]{Chenfei Wang $^{\ast\ast}$}\ead{Wangcf25@mails.jlu.edu.cn}              

\address[Paestum]{School of Mathematics, Jilin University,
130012, Changchun, China   }  


\begin{abstract}                          
This paper investigates discrete-time stochastic recursive games driven by fractional noises, in which the cost functional is defined through a backward stochastic difference equation. By introducing discrete-time adjoint equations, we establish a maximum principle and the verification theorem for open-loop Nash equilibria in non-zero sum and zero sum stochastic recursive game problems. The effectiveness of our theorem is demonstrated through applications to two different linear-quadratic stochastic recursive games and a portfolio selection problem with a recursive utility functional.
\end{abstract}
\end{frontmatter}

\section{Introduction}

Stochastic differential games have long served as a foundational framework for modeling strategic interactions in dynamic and uncertain environments. Originating from the seminal work of Isaacs \cite{IsaacsA} on deterministic differential games. In continuous-time setting, foundational contributions by Fleming and Souganidis \cite{WOnthexistencegames} established the existence of value functions and their characterization via viscosity solutions of Hamilton–Jacobi–Bellman–Isaacs (HJBI) equations. Subsequent extensions incorporated backward stochastic differential equations (BSDEs) by Pardoux and Peng \cite{PARDOUX1990ADAPTEDDIFFERENTIAL-EQUATION}, which facilitated a more natural formulation of dynamic optimization under stochastic terminal conditions and enriched applications in finance, economics, and insurance. The application of BSDE methods to game theory, pioneered by Hamadène and Lepeltier \cite{1995Zero} among others, further simplified the analysis of continuous-time stochastic differential games. Further extensions and applications of continuous-time stochastic differential games have  been explored. Wang and Yu \cite{Wangapplications} investigated a kind of non-zero sum differential game of BSDEs. Huang et al. \cite{Huangnoise.} studied a Linear-quadratic (LQ) mean-field game of stochastic large-population system with partial information and common noise. Lin  \cite{Linbarriers} investigated Nash equilibrium payoffs for nonzero-sum stochastic differential games with reflection. Li and Li \cite{LiCONDITION}  studied on the case of zero-sum games with asymmetric information on the pay-off and the case of symmetric information non-zero sum differential game. 
\par These above results are mainly concentrated on the continuous-time problems. The discrete-time problem is as important as the continuous-time one. 
Discrete-time stochastic games are frequently investigated within the frameworks of Markov decision processes through stochastic programming (see, e.g., Başar and Olsder \cite{Başaynamicnoncooperative}). This approach has established a theoretical basis for modeling sequential decision-making under uncertainty. Pachter and Pham \cite{MDynamic} considered the application of the theory of discrete-time LQ dynamic games to digital signal processing. Mukaidani and Xu \cite{Mukaidanisystems} studied infinite horizon LQ Stackelberg games for a discrete-time stochastic system with multiple decision makers. Sun et al. \cite{Sundiscrete-time} dealed with the infinite horizon LQ differential games for discrete-time stochastic systems with both state and control dependent noises. While the mentioned results are established within the LQ framework, this paper extends the investigation to a general setting for the discrete-time stochastic game problems.

The stochastic recursive control problem, introduced by Pardoux and Peng \cite{PARDOUX1990ADAPTEDDIFFERENTIAL-EQUATION}, is characterized by a state governed by a controlled stochastic differential equation (SDE) and a cost functional defined by the solution of a coupled controlled BSDE. It generalizes the classical stochastic control problem, which appears as the special case where the BSDE's generator is independent of its solution pair. This framework, equivalent to the stochastic differential utility (SDU) proposed by Duffie and Epstein \cite{Darrell1992Stochastic}, naturally models dynamic recursive preferences and has made the optimal consumption-investment problem a central topic in mathematical finance.

While recursive structures have been extensively studied in stochastic control theory, they have received considerably less attention in the context of stochastic differential games. Motivated by this gap, this paper focuses on discrete-time stochastic recursive game problems.

The majority of existing studies on stochastic games concentrate on systems driven by standard Brownian motion. In this setting, the dynamic programming principle and the associated Hamilton\allowbreak –Jacobi\allowbreak –Bellman\allowbreak –Isaacs equations provide a rigorous framework for deriving optimal feedback controls, as developed in the foundational works on stochastic optimal control by Fleming and Rishel \cite{Fleming1975Deterministic} and Yong and Zhou \cite{1999Stochastic}, and further extended to the game-theoretic setting by Wu and Zhang \cite{WuGAME}. In contrast, the non-Markovian property of fractional noise renders the dynamic programming approach inapplicable. This necessitates the adoption of the maximum principle as the fundamental solution tool. Bayraktar and Poor \cite{Bayraktar2005Stochastic} investigated stochastic differential games under non-Markovian settings. However, it has primarily addressed linear cases as illustrative examples. Therefore, our work fills this gap by establishing a discrete-time maximum principle for stochastic recursive game problems with fractional noise under a general setting.

Another significant highlight is as follows. In the continuous-time setting, Han et al. derived necessary optimality conditions for controlled systems driven by fractional Brownian motion (see Theorem 5.7 of \cite{YuecaiH}). However, their conditions involve two equalities, one of which is $\sigma_u^*(t)p_t=0$. This  makes their framework almost inapplicable to the cases in which the diffusion coefficient depends nontrivially on the control variables. In contrast, in this work, we derive a single equality as the necessary optimality condition for each control process. Furthermore, the optimal control obtained in the examples explicitly depends on past information of the fractional noise. This illustrates the practical relevance of our theoretical results.

The paper is structured as follows. In Section 2, we introduce some notions and formulation of the stochastic game problems. In section 3, we derive the maximum principles for two distinct types of discrete-time stochastic recursive game problems: the discrete-time non-zero sum stochastic recursive game problem (\textbf{Problem (NZSRG)}) and the discrete-time zero sum stochastic recursive game problem (\textbf{Problem (ZSRG)}). Finally, in section 4, three concrete examples are employed to validate the proposed stochastic maximum principle, specifically, two types of LQ stochastic game problems (\textbf{Problem (LQNZSRG)} and \textbf{Problem (LQZSRG)}) and a portfolio selection problem with a recursive utility function.   

\section{Preliminary}
\subsection{Basic Notions}
\par Let $H\in(0,1)$ be a fixed constant. 
   $ B^H=(B^H_1(t),\cdots,\\B^H_m(t),t\geq0)$
 be a $m$-dimensional fractional Brownian motion of Hurst parameter $H$. The fractional noise $\xi_n^H=B^H(n+1)-B^H(n),(n=0,\cdots,N-1)$ is a mean 0 Gaussian process with the following covariance 
\begin{equation}\label{1}
    \begin{split}
&\rho_H(i,j)
:=\mathbf{E}\big[\xi_i^H\xi_j^H\big]\\
        =&\dfrac{1}{2}\bigg[\big|i-j+1\big|^{2H}+\big|i-j-1\big|^{2H}-2\big|i-j\big|^{2H}\bigg].
    \end{split}
\end{equation}

\par Let the filtration $\mathbb{F}=(\mathcal{F}_n)_{0\leq n\leq N}$ defined by $\mathcal{F}_n=\mathcal{F}_0\vee \sigma(\xi_0,\xi_1,\cdots,\xi_{n-1}),$ and $(\Omega,\mathcal{F},\{\mathcal{F}_n\}_{0\leq n\leq N},\mathbb{P})$ be a filtered probability space, $\mathcal{F}_0\subset\mathcal{F}$ be a sub $\sigma-$algebra. 
$x$ is a given $\mathcal{F}_0$-measurable random variables taking values in $\mathbf{R}^n$. We denote $\mathcal{T}=\{0,1,\dots,N-1\}$.

   For $\beta\geq 1$, denote the set of $\mathbf{F}$-adapted  process $X=(X_n)_{0\leq n\leq N}\in L_{\mathbf{F}}^{\beta}(0,N)$ such that 
\begin{equation*}
\big\|X\big\|_{\beta}=\bigg(\sum_{n=0}^N\mathbf{E}\big\|X_n\big\|^{\beta}\bigg)^{\tfrac{1}{\beta}}<\infty.
\end{equation*}
According to the Han and Li  \cite{Yuecai}, we have the following lemma.
\begin{lem}\label{9}
Let the matrix $E = \big(e(i,j)\big)_{i,j\geq 0}$ be defined recursively by the formula
\begin{equation}\label{2}
e(i,j)= \mathbf{1}_{\{j \leq i\}} \cdot \frac{\rho_H(i,j) - \sum_{m=0}^{j-2} e(i,m+1)e(j,m+1)}{e(j,j)},\end{equation}
and we define the coefficients $c(\cdot, \cdot)$ via $c(i-1,j-1) = (E^{-1})_{ij}$.
    Then the sequence $\{\eta_n\}$ defined by
\begin{equation*}
       \eta_n = \sum_{m=0}^{n} c(n,m) \, \xi_m^H, \quad n \geq 0, 
\end{equation*}
forms a sequence of independent Gaussian random variables such that each $\eta_n$ is $\mathcal{F}_{n+1}$-measurable and independent of $\mathcal{F}_n$.  
We refer to $\{\eta_n\}$ as the ``underlying white noise'' of the process $\{\xi_n\}$. Moreover, we have
\begin{equation*}
\quad \xi_n^H = \sum_{m=0}^{n} e(n,m) \, \eta_m.
\end{equation*}
\end{lem}

\subsection{Stochastic Recursive Game formulations}\label{24}

 The system is described by 
\begin{equation}\label{3}
    \begin{split}
    \begin{cases}
X(n+1)=X(n)+b\big(n,X(n),u_{1}(n),u_{2}(n)\big)\\
\quad\quad\quad\quad\quad+\sigma\big(n,X(n),u_{1}(n),u_{2}(n)\big)\xi_n^H,\\
X(0)=x,\quad\quad\quad\quad\quad\quad\quad\quad n\in\mathcal{T},
    \end{cases}
    \end{split}
\end{equation}
where $x\in\mathbf{R}^d$ is $\mathcal{F}_0$-measurable and $u_i,i=1,2$ is the control processes. $b(n,\cdot,\cdot,\cdot):\mathcal{T}\times\mathbf{R}^d\times \mathbf{R}^{k_1}\times \mathbf{R}^{k_2}\rightarrow\mathbf{R}^{d},\sigma(n,\cdot,\cdot,\cdot):\mathcal{T}\times\mathbf{R}^d\times \mathbf{R}^{k_1}\times \mathbf{R}^{k_2}\rightarrow\mathbf{R}^{d\times m}.$
 
 Let $U_i \subset \mathbf{R}^{k_i}$ be a nonempty convex set. The open-loop admissible control set for Player $i$ is defined as
    \begin{equation*}
    \begin{split}
    \mathbf{U}_i := \big\{ u_i \in L_{\mathcal{F}}^{2a}(0, N-1; 
    \mathbf{R}^{k_i}) \;|\; u_i(n) \in U_i, \\   
          \text{ a.s. for all }
         n \in \mathcal{T} \big\}, \quad i=1,2.
    \end{split}
    \end{equation*}
    \begin{rem}
        Let $a>1$ is a fixed constant. In this paper, we frequently use the spaces $L^{2a}$
  with $a
>
1$. The reason for requiring this condition is that fractional noises exhibit
autocorrelation for all \(H\in(0,1)\); in particular, for \(H>1/2\), they
exhibit long-range dependence $($see Proposition 1 of Han and Li \cite{Yuecai}$)$. To handle the correlations in the noise, we need higher-order moment estimates beyond the usual $L
^2$-framework. The higher integrability provided by $L 
^{2a}$
  allows us to apply Hölder's inequality and an interpolation argument $($with a parameter $\delta\in[a^{-1/N},1)$$)$ to derive the desired $L^ 
2$-bounds for the state, adjoint, and variation processes. 
    \end{rem}
    
 Let $u^{\ast}_i$ be the optimal control and $\bar{u}_i$ be any other admissible controls. Denote $v_i=\bar{u}_i-u_i^{\ast}$. For any $\theta\in(0,1)$, we define $u^{\theta}_i=u^{\ast}_i+\theta v_i$. 

    Denote by $J_i(u_1, u_2)$ the cost functional of Player $i$. For simplicity, in this subsection, we only consider $i=2$ the case of two players, which is similar for $n$ players. 

  \textbf{Problem (NZSRG)}  This problem is to find $(u_1^{\ast}, u_2^{\ast}) \in \mathbf{U}_1\times\mathbf{U}_2$ that satisfies the following equilibrium conditions,
\begin{equation}\label{31}
        \begin{split}
            J_1(u_1^{\ast}, u_2^{\ast}) = \inf_{u_1 \in \mathbf{U}_1} J_1(u_1, u_2^{\ast}), \\
            J_2(u_1^{\ast}, u_2^{\ast}) = \inf_{u_2 \in \mathbf{U}_2} J_2(u_1^{\ast}, u_2),
        \end{split}
    \end{equation}
and let the cost functional satisfy
\begin{equation*}
    J_i(u_i, u_{-i}) = Y_i(0), \quad i = 1,2.
\end{equation*}
Here $\big(Y_i(n), Z_i(n)\big)$ is the solution to the following BS$\Delta$E,
\begin{equation}\label{eq:adjoint}
\begin{cases}
Y_i(n)+Z_i(n)\eta_n
=Y_i(n+1)+f_i\big(n+1,\psi_i(n+1)\big),\\[4pt]
Y_i(N)=g_i\big(X(N)\big),
\end{cases}
\end{equation}
where $\psi_i(n+1):=\big(X^i(n+1),Y_i(n+1),Z_i(n+1),u_1(n+1),u_2(n+1)\big),X^i(n)=X^i\big(n,X^i(n),u_i(n),u_{-i}^{\ast}(n)\big)$ 
and $f_i(n,\cdot,\cdot,\cdot,\cdot,\cdot):\mathcal{T}\times\mathbf{R}^d \times\mathbf{R}^{n}\times\mathbf{R}^{n\times m}\times \mathbf{R}^{k_1}\times \mathbf{R}^{k_2}\rightarrow\mathbf{R}^{n}$, $g_i(\cdot):\mathbf{R}^d\rightarrow \mathbf{R}^{n}$. Such a pair $(u_1^*, u_2^*)$ is called an open-loop Nash equilibrium point, and the corresponding state process is denoted by $X^* = \{X^*(n),n\in\mathbf{T}\}$.

\textbf{Problem (ZSRG)} Under the state equation (\ref{3}), the problem here is to find $(u_1^{\ast},u_2^{\ast})\in \mathbf{U}_1\times \mathbf{U}_2$ satisfying
\begin{equation}\label{17}
    J(u_1^{\ast},u_2)\leq J(u_1^{\ast},u_2^{\ast})\leq J(u_1,u_2^{\ast}),\quad (u_1,u_2)\in \mathbf{U}_1\times \mathbf{U}_2,
\end{equation}
where 
\begin{equation}\label{y_0}
J(u_1,u_2)=Y(0),
\end{equation}
and $\big(Y(n),Z(n)\big)$ is the solution to the following BS$\Delta$E,
\begin{equation}\label{26}
    \begin{split}
        \begin{cases}
            Y(n)+Z(n)\eta_n=Y(n+1)+f\big(n+1,\psi(n+1)\big),\\
            Y(N)=g\big(X(N)\big).
        \end{cases}
    \end{split}
\end{equation}

The pair of Nash equilibrium point $(u_1^{\ast},u_2^{\ast})$ satisfying (\ref{17}) is called a saddle point, and the corresponding state is denoted by $X^{\ast}=\big\{X^{\ast}(n),n\in[0,N]\big\}$.
    In order not to over complicate the already notational heavy presentation of this paper, in what follows we shall assume all processes are 1-dimensional (i.e., $d=n=k_1=k_2=m=1$).

\section{Maximum principle for two kinds of discrete-time stochastic games }
According to the state equations and recursive cost functionals (\ref{3})--(\ref{26}) of section \ref{24}. Assume that the coefficients $b,\sigma,f,g,f_i,g_i $ $(i=1,2)$ satisfy the assumptions (H1)--(H3) with $\mu$ replaced by $(u_1,u_2)$, and for all $(n,x,y,\mu)\in\mathcal{T}\times\mathbf{R}\times\mathbf{R}\times\mathbf{R}^2.$
\par (H1) For all $(n,x,\mu)$, $b(\cdot,x,\mu), \sigma(\cdot,x,\mu)\in L_{\mathcal{F}}^{2a}(0,N-1)$, and there exists $C>0$,
\begin{equation*}
    \begin{split}
        &\left|\varphi(n,x_1,\mu_1)-\varphi(n,x_2,\mu_2)\right|
        \leq C\big(\big|x_1-x_2\big|+\big|\mu_1-\mu_2\big|\big),\\
    \end{split}
\end{equation*}
for $\varphi=b,\sigma,b_x,\sigma_x,b_\mu,\sigma_\mu.$

\par (H2) For all $(x,y,z,\mu)$, $f,f_i(\cdot,x,y,z,\mu)\in L_{\mathcal{F}}^{2a}(0,N-1)$, $g,g_i(x)\in L^{2a}(\mathcal{F}_N)$. There exists $L>0$, such that 
\begin{equation*}
    \begin{split}
        &|\varphi^f_{x}(n,x_1,y_1,z_1,\mu_1)-\varphi^f_{x}(n,x_2,y_2,z_2,\mu_2)|\\
        \leq& L\bigg(\big|x_1-x_2\big|+\big|y_1-y_2\big|+\big|z_1-z_2\big|+\big|\mu_1-\mu_2\big|\bigg),\\
        &|\varphi^g_{x}(x_1)-\varphi^g_{x}(x_2)|\le L|x_1-x_2|,
        \\&\varphi^f=f,f_i,\quad\varphi^g=g,g_i.
    \end{split}
\end{equation*}

\par  (H3) $b(N,x,\mu)=\sigma(N,x,\mu)=f(0,x,\mu)=f_i(0,x,\mu)=0,$ for all $(x,\mu)$.

(H4) For all $(n,x,y,\mu)$, the functions $b(n,\cdot,\cdot)$, $\sigma(n,\cdot,\cdot)$, $f(n,\cdot,\cdot,\cdot,\cdot),f_i(n,\cdot,\cdot,\cdot,\cdot),g(\cdot)$ and $g_i(\cdot)$ are continuously differentiable with respect to $(x,\mu)$, and their partial derivatives satisfy 
\begin{equation*}
\begin{split}
&\big|b_x(n,x,\mu)\big|+\big|b_v(n,x,\mu)\big|
\\+&\big|\sigma_x(n,x,\mu)\big|+\big|\sigma_v(n,x,\mu)\big|\le C,\\
&\big|\varphi^f_{x}(n,x,\mu)\big|+\big|\varphi^f_{\mu}(n,x,\mu)\big|
+\big|\varphi^g_{x}(x)\big|\\\le& C\big(1+|x|+|\mu|\big),\\&\varphi^f=f,f_i,\quad\varphi^g=g,g_i,
\end{split}
\end{equation*}
for some constant $C>0$.
\subsection{Maximum principle for discrete-time non-zero sum stochastic game}
We introduce the variational equation as follows, 
\begin{equation}\label{25}
\begin{cases}
V^i(n+1) =V^i(n)+ \left[ b_{x}^{*}(n) V^i(n) + b_{u_i}^{*}(n) v_{i}(n) \right] \\\quad\quad\quad\quad\quad+ \left[ \sigma_{x}^{*}(n) V_i(n) + \sigma_{u_i}^{*}(n) v_{i}(n) \right] \xi_n^H, \\
V^i(0) = 0, \quad i=1,2,
\end{cases}
\end{equation}
where
$\phi^i(n)=\phi^i\big(n,X(n),u_i(n),u_{-i}^{\ast}(n)\big),\phi^{*}(n)=\phi\big(n,X^{\ast}(n),u_i^{\ast}(n),u_{-i}^{\ast}(n)\big),\quad\phi=b,\sigma,b_x,b_{u_i},\sigma_x,\sigma_{u_i}.$

By Lemma 2 of Han and Li \cite{Yuecai}, S$\Delta$E (\ref{25}) has a unique solution $V^i\in L_{\mathcal{F}}^{2}(0,N)$. 

Furthermore, consider the associated linear BS$\Delta$E,
\begin{equation}\label{4}
    \begin{split}
        \begin{cases}
            R_{i}(n)+S_{i}(n)\eta_n=R_{i}(n+1)+f_{i,x}^{\ast}(n+1)V^{i}(n+1)\\
            \quad\quad\quad+f_{i,y}^{\ast}(n+1)R_{i}(n+1)+f_{i,u_i}^{\ast}(n+1)v_{i}(n+1)\\\quad\quad\quad+f_{i,z}^{\ast}(n+1)S_i(n+1),\quad\quad n\in\mathcal{T}\\
            R_i(N)=g_{i,x}\big(X(N)\big)V^i(N), \quad\quad i=1,2,
        \end{cases}
    \end{split}
\end{equation}
where $f_i^{\ast}(n)=f_i\big(n,X^{\ast}(n),Y_i^{\ast}(n),Z_i^{\ast}(n),u_i^{\ast}(n),u_{-i}^{\ast}(n)\big).$
\begin{rem}
    In this paper, we use the following notation,
\begin{equation*}
    (u_i,u_{-i}^{\ast}):=\bigg\{\begin{array}{cc}
       (u_1,u_2^{\ast}),\quad i=1,\\
(u_1^{\ast},u_2),\quad i=2.
    \end{array}
\end{equation*}
Moreover, the superscript \(i\) attached to the functions and processes is used to indicate a unilateral deviation of the \(i\)-th player, with \(u_{-i}^{\ast}(n)\) being the optimal strategy of the other player in the two-player game held fixed. Correspondingly, the superscript \(\ast\) on the left-hand side means that these functions and processes are evaluated at the fully optimal variables \((X^{\ast},u_i^{\ast},u_{-i}^{\ast})\).
\end{rem}

By Theorem 2.6 of Ji and Zhang \cite{Yuanyuan}, BS$\Delta$E (\ref{4}) has a unique solution $(R_i,S_i)\in L_{\mathcal{F}}^{2}(0,N)\times L_{\mathcal{F}}^{2}(0,N-1) $.
\par Based on Lemma 3 of Han and Li \cite{Yuecai}, we establish the necessary estimates as follows.
\begin{lem}\label{10}
    Let $X^{\theta},X^{\ast}$ be the corresponding state equation to $u_i^{\theta},u_i^{\ast}$, respectively. According to the Assumptions (H1)--(H4), we have 
\begin{equation}\label{6}
\sum_{n=0}^N\mathbf{E}\big|X^{i,\theta}(n)-X^{*}(n)\big|^2=O(\theta^2),\quad i=1,2,
  \end{equation}
     and as $\theta\rightarrow 0,$ 
\begin{equation}\label{7}
        \lim_{\theta\rightarrow0}\sum_{n=0}^N\mathbf{E}\bigg|\dfrac{X^{i,\theta}(n)-X^{*}(n)}{\theta}-V^i(n)\bigg|^2=0,\quad i=1,2,
    \end{equation}
   where
    \begin{equation*}
        \begin{split}
        \begin{cases}
              X^{i,\theta}(n+1)=X^{i,\theta}(n)+b^i\big(n,X^{\theta}(n),u_{i}^{\theta}(n),u_{-i}^{\ast}(n)\big)\\
\quad\quad\quad\quad\quad+\sigma^i\big(n,X^{\theta}(n),u_{i}^{\theta}(n),u_{-i}^{\ast}(n)\big)\xi_n^H,\\X^{*}(n+1)=X^{*}(n)+b^i\big(n,X^{\ast}(n),u_{i}^{\ast}(n),u_{-i}^{\ast}(n)\big)\\
\quad\quad\quad\quad\quad+\sigma^i\big(n,X^{\ast}(n),u_{i}^{\ast}(n),u_{-i}^{\ast}(n)\big)\xi_n^H,
n\in\mathcal{T}.
\end{cases}
        \end{split}
    \end{equation*} 
\end{lem}

$Proof.$ The proof of equality (\ref{6}) is based on Lemma 5 from the work of Han and Li \cite{Yuecai}.
\par For (\ref{7}), let $\hat{V}^{i,\theta}(n)=\dfrac{X^{i,\theta}(n)-X^{*}(n)}{\theta}-V^i(n)$, 
then for each $n\in\mathcal{T}$ and each $\theta>0$, we have
    \begin{equation*}
        \begin{split}
&\hat{V}^{i,\theta}(n+1)\\
        =&\hat{V}^{i,\theta}(n)+ \bigg\{b_{x}^{\lambda,*}(n)\hat{V}i^{i,\theta}(n)+\Big[b_{x}^{i,\lambda,\theta}(n)-b_{x}^{*}(n)\Big]\\
        &\cdot V^{i}(n)+\Big[b_{u_i}^{i,\lambda,\theta}(n)-b_{u_i}^{*}(n)\Big]v_i(n)\bigg\}\\
        & +\bigg\{\sigma_{x}^{\lambda,*}(n)\hat{V}^{i,\theta}(n)+\bigg[\sigma_{x}^{i,\lambda,\theta}(n)-\sigma_{x}^{*}(n)\bigg]V^i(n)\\&+\bigg[\sigma_{u_i}^{i,\lambda,\theta}(n)-\sigma_{u_i}^{*}(n)\bigg]v_{i}(n)\bigg\}\xi_n^H,\\
        \end{split}
    \end{equation*}
 where 
 \begin{equation*}
     \begin{split}
       \phi^{\lambda,\theta}(n)=\int_0^1\phi\big(n,X^{\ast}(n)+\lambda(X^{\theta}(n)-X^{\ast}(n)),\\u_i^{\ast}(n)+\lambda\theta v_i(n),u_{-i}^{\ast}(n)\big)d\lambda,\\\phi=b,\sigma,b^i_{x},\sigma^i_{x},b^i_{u_i},\sigma^i_{u_i}.
        \end{split}
 \end{equation*}
 Then  
\begin{equation*}
     \begin{split}
          &\mathbf{E} \left| \hat{V}^{i,\theta}(n+1)\right|^{2a\delta^{n+1}}
\\\leq& C \mathbf{E} \left| \hat{V}^{i,\theta}(n) \right|^{2a\delta^{n+1}}
+ C \mathbf{E} \left\| b_{u_i}^{i,\lambda,\theta}(n) - b_{u_i}^{*}(n) \right\|^{2a\delta^{n+1}}\\&\cdot \left| v_i(n) \right|^{2a\delta^{n+1}} + C \mathbf{E} \left\| b_{x}^{i,\lambda,\theta}(n) - b_{x}^{*}(n) \right\|^{2a\delta^{n+1}} \notag \\
& \cdot\left| \theta^{-1}\big(X^{i,\theta}(n) - X^{*}(n)\big) \right|^{2a\delta^{n+1}}+ C \mathbf{E}\notag \\
& \cdot\bigg[ \left| \hat{V}^{i,\theta}(n) \right|^{2a\delta^{n+1}}+
  \left\| \sigma_{u_i}^{i,\lambda,\theta}(n) - \sigma_{u_i}^{*}(n) \right\|^{2a\delta^{n+1}} \\&\cdot\left| v_i(n) \right|^{2a\delta^{n+1}} + \left\| \sigma_{x}^{i,\lambda,\theta}(n) - \sigma_{x}^{*}(n) \right\|^{2a\delta^{n+1}}  \notag \\
& \cdot\left| \theta^{-1}\big(X^{i,\theta}(n) - X^{*}(n)\big) \right|^{2a\delta^{n+1}} \bigg] \left| \xi_{n}^{H} \right|^{2a\delta^{n+1}} \notag \\
\leq& C \left[ \left( \mathbf{E} \left| \hat{V}^{i,\theta} (n)\right|^{2a\delta^{n}} \right)^{\delta} + o(1) \right]
 \notag \\
& \cdot\left[ 1 + \left( \mathbf{E} \left| \xi_{n}^{H} \right|^{\frac{2a\delta^{n+1}}{1-\delta}} \right)^{1-\delta} \right],
\end{split}
 \end{equation*}
 where $\delta\in(\frac{1}{a},1)$.
Since   $\hat{V}^{i,\theta}(0)=0, \, a.s.,$
it is a straightforward result that
\begin{equation*}
    \sum_{n=0}^N\mathbf{E}\big|\hat{V}^i(n+1)\big|^2\leq \sum_{n=0}^N\mathbf{E}\big(\big|\hat{V}^i(n+1)\big|^{2a\delta^{n}}\big)^{\tfrac{1}{a\delta^n}}=o(1).
\end{equation*}
This completes the proof of (\ref{7}).
\begin{lem}\label{46}
Under Assumptions (H1)--(H4), the following estimates hold for the solution pairs
\((Y_i^{\theta}, Z_i^{\theta})\) and \((Y_i^{*}, Z_i^{*})\)
corresponding to the controls \(u_i^{\theta}\) and \(u_i^{*}\), respectively, \begin{equation}\label{49}
        \sum_{n=0}^N \mathbf{E} \left[ \left| Y_i^{\theta}(n) - Y_i^{*}(n) \right|^2 + \left| Z_i^{\theta}(n) - Z_i^{*}(n) \right|^2 \right] = O(\theta^2).
    \end{equation}
Furthermore, defining the remainder term of the variation
\begin{equation*}
        \begin{split}
            \hat{Y}_i^{\theta}(n) = \frac{Y_i^{\theta}(n) - Y_i^{*}(n)}{\theta} - R_i(n), \\
        \hat{Z}_i^{\theta}(n) = \frac{Z_i^{\theta}(n) - Z_i^{*}(n)}{\theta} - S_i(n),
        \end{split}
    \end{equation*}
    then we have
\begin{equation}\label{50}
        \lim_{\theta \to 0} \sum_{n=0}^N \mathbf{E} \left[ \left| \hat{Y}_i^{\theta}(n) \right|^2 + \left| \hat{Z}_i^{\theta}(n) \right|^2 \right] = 0.
    \end{equation}
\end{lem}

$Proof.$  The proof of the equality (\ref{49}) proceeds by backward induction, applying Lemma \ref{10}. Under the inductive hypothesis that
\begin{equation*}
    \begin{split}
        &\mathbf{E}\big|Y_{i}^{\theta}(n)-Y_{i}^{\ast}(n)\big|^2+\mathbf{E}\big|Z_{i}^{\theta}(n)-Z_{i}^{\ast}(n)\big|^2\big|\eta_n\big|^2\\
        \leq 3&C\mathbf{E}\bigg[\big|Y_{i}^{\theta}(n+1)-Y_{i}^{\ast}(n+1)\big|^2\\+&\big|Z_{i}^{\theta}(n+1)-Z_{i}^{\ast}(n+1)\big|^2+\big|u_{i}^{\theta}(n+1)-u_{i}^{\ast}(n+1)\big|^2\\+&\big|X^{i,\theta}(n+1)-X^{*}(n+1)\big|^2\bigg]=O(\theta^2).
    \end{split}
\end{equation*}

Then we use the same method of Lemma \ref{10}. 
Therefore, it turns out that
\begin{equation*}
    \lim_{\theta\rightarrow 0}\mathbf{E}\big|\hat{Y}_i^{\theta}(n)\big|^{2}+\mathbf{E}\big|\hat{Z}_i^{\theta}(n)\big|^{2}=0,
\end{equation*}
for $n\in\mathcal{T}$.  Then we conclude 
\begin{equation*}
    \lim_{\theta\rightarrow0}\sum_{n=0}^N\mathbf{E}\bigg[\big|\hat{Y}_i^{\theta}(n)\big|^2+\big|\hat{Z}_i^{\theta}(n)\big|^2\bigg]=0,
\end{equation*}
which completes the proof. 
\par Next, we introduce the adjoint S$\Delta$E as follows,
\begin{equation}\label{51}
    \begin{split}
        \begin{cases}
            k_i(n+1)=k_i(n)+f_{i,y}^{\ast}(n)k_i(n)+f_{i,z}^{\ast}(n)k_i(n)\eta_n,\\
        k_{i}(1)=-1, k_{i}(0)=0.
        \end{cases}
    \end{split}
\end{equation}
Obviously, the equation (\ref{51})  has a unique solution $k_i\in L_{\mathcal{F}}^2(0,N)$.
\par Now, we have the following two results for the equilibrium point of \textbf{Problem (NZSRG)}. 
\begin{thm}\label{19}
    Suppose $(u_1^{\ast},u_2^{\ast})$ is an equilibrium point of \textbf{Problem (NZSRG)}, $X^{\ast}=\big\{X^{\ast}(n),n\in[0,N] \big\}$ is the corresponding state, and the following adjoint equation is the dual equation of S$\Delta$E (\ref{25}), i=1,2. 
    \begin{equation}\label{8}
        \begin{split}
            \begin{cases}
        p_{i}(n)+q_{i}(n)\eta_n=p_{i}(n+1)+b_{x}^{*}(n+1)p_{i}(n+1)\\\quad\quad\quad\quad\quad\quad\quad+e(n+1,n+1)\sigma^{*}_{x}(n+1)q_{i}(n+1)\\
        \quad\quad\quad\quad\quad\quad\quad+\sigma_{x}^{*}(n+1)p_i(n+1)\mathbf{E}^{\mathcal{F}_{n+1}}\big[\xi_{n+1}^H\big]\\\quad\quad\quad\quad\quad\quad\quad-f_{i,x}^{\ast}(n+1)k_{i}(n+1),n\in\mathcal{T}\\
            p_{i}(N)=-g_{i,x}(X^{\ast}(N))k_i(N).
            \end{cases}
        \end{split}
    \end{equation}
    Then the following inequality holds, 
\begin{equation}\label{inequa1}
       \begin{split}
           &\Big[b_{u_1}^{*}(n)p_{1}(n)+\sigma_{u_1}^{*}(n)p_{1}(n)\sum_{m=0}^{n-1}d(n,m)\xi_m^H\\&+e(n,n)\sigma_{u_1}^{*}(n)q_{1}(n)-f_{1,u_1}^{\ast}(n)k_{1}(n)\Big]\\
\cdot&\big(u_1(n)-u_1^{\ast}(n)\big)\geq 0,\quad \text{ a.s.,\, for all $n\in\mathcal{T}, u_i\in\mathbf{U}_i$},\end{split}\end{equation}
where $d(n,m)=\sum_{l=0}^{n-1}e(n,l)c(l,m).$
    \end{thm} 
    
  $Proof.$  According to the assumptions (H1)--(H4) and the Theorem 3 in the paper of Han and Li \cite{Yuecai}, it is easy to check  BS$\Delta$E (\ref{8}) has a unique solution. 
Taking summation and expectation, we have
\begin{equation}\label{13}
    \begin{split}
        &\mathbf{E}\sum_{n=0}^{N-1}\Delta(p_1(n)V^1(n))+\mathbf{E}\sum_{n=1}^{N-1}\Delta(R_1(n)k_1(n))\\
=&\mathbf{E}\bigg\{\sum_{n=0}^{N-1}\sigma_{u_1}^{*}(n)q_1(n)v_1(n)\eta_n\xi_n^H\\&+\sum_{n=0}^{N-1}\sigma_{u_1}^{*}(n)p_1(n)v_1(n)\xi_n^H\\&+\sum_{n=1}^{N}f_{1,x}^{\ast}(n)V^1(n)k_1(n)+\sum_{n=0}^{N-1}b_{u_1}^{*}(n)p_1(n)v_1(n)\bigg\}\\
    &+\mathbf{E}\bigg\{-\sum_{n=1}^Nk_1(n)\big[f_{1,x}^{\ast}(n)V^1(n)+f_{1,u_1}^{\ast}(n)v_1(n)\big]\\&+\big[\big(R_1(0)-R_1(1)\big)k_1(1)\big]\bigg\},
    \end{split}
\end{equation}
where  $b_{x}^{*}(N)=\sigma_{x}^{*}(N)=f_{1,x}^{\ast}(N)=V^1(0)=0$ and $f_{1,y}^{\ast}(N)=f_{1,z}^{\ast}(N)=k_1(0)=0.$

According to the Lemma 1, we have
\[
\xi_n^H=\sum_{m=0}^{n}e(n,m)\eta_m,
\]
where \(\eta_m\) is \(\mathcal{F}_{m+1}\)-measurable and independent of \(\mathcal{F}_m\). Taking the conditional expectation yields
\[
\mathbf{E}^{\mathcal{F}_n}[\xi_n^H]
=\sum_{m=0}^{n-1}e(n,m)\eta_m,
\]
and there exist coefficients \(d(n,m)\) such that
\[
\sum_{m=0}^{n-1}e(n,m)\eta_m
=\sum_{m=0}^{n-1}d(n,m)\xi_m^H,
\]
where $d(n,m)=\sum_{l=0}^{n-1}e(n,l)c(l,m)$, $c(\cdot,\cdot),e(\cdot,\cdot)$ are given by Lemma \ref{9}.
Therefore,
\[
\mathbf{E}^{\mathcal{F}_n}[\xi_n^H]
=\sum_{m=0}^{n-1}d(n,m)\xi_m^H.
\]
Then we have
\begin{equation}\label{23}
    \begin{split}&\mathbf{E}\big[\sigma_{u_1}^{*}(n+1)p_1(n+1)v_1(n+1)\xi_n^H\big]\\
=&\mathbf{E}\bigg[\sigma_{u_1}^{*}(n)p_1(n)v_1(n)\sum_{m=0}^{n-1}d(n,m)\xi_m^H\bigg],
 \end{split}
\end{equation}
and 
\begin{equation}\label{12}
    \begin{split}
&\mathbf{E}\big[\sigma_{u_1}^{*}(n)q_1(n)v_1(n)\eta_n\xi_n^H\big]\\
=&\mathbf{E}\big[e(n,n)\sigma^{*}_{u_1}(n)q_1(n)v_1(n)\big].
    \end{split}
\end{equation}
Then we get 
\begin{equation*}
    \begin{split}
        k_1(1)R_1(1)+\mathbf{E}\bigg\{\sum_{n=0}^{N-1}v_1(n)\Big[e(n,n)\sigma^{*}_{u_1}(n)q_1(n)&\\+b^{*}_{u_1}(n)p_1(n)+\sigma_{u_1}^{*}(n)p_1(n)\sum_{m=0}^{n-1}d(n,m)\xi_m^H&\\-f_{1,u_1}^{\ast}(n)k_1(n)\Big]+\big[\big(R_1(0)-R_1(1)\big)k_1(1)\big]\bigg\}&=0.
    \end{split}
\end{equation*}Therefore, 
\begin{equation*}
    \begin{split}
        &\frac{d}{d\theta}J_1\big(u_1^{\theta},u_2^{\ast}\big)\\=&R_1(0)\\=&\mathbf{E}\bigg\{\sum_{n=0}^{N-1}\Big[e(n,n)\sigma^{*}_{u_1}(n)q_1(n)\\&+\sigma_{u_1}^{*}(n)p_1(n)\sum_{m=0}^{n-1}d(n,m)\xi_m^H\\&+b_{u_1}^{*}(n)p_1(n)-f_{1,u_1}^{\ast}(n)k_1(n)\Big]v_1(n)\bigg\}.
    \end{split}
\end{equation*}
By the arbitrary of $v_1$, we get 
\begin{equation*}
\begin{split}
      &\mathbf{E}\bigg\{\textbf{1}_{\mathcal{A}}\Big[b_{u_1}^{*}(n)p_{1}(n)+\sigma_{u_1}^{*}(n)p_{1}(n)\sum_{m=0}^{n-1}d(n,m)\xi_m^H\\&+e(n,n)\sigma_{u_1}^{*}(n)q_{1}(n)-f_{1,u_1}^{\ast}(n)k_{1}(n)\Big]\\
      \cdot&\big(u_1(n)-u_1^{\ast}(n)\big)\bigg\}\geq 0,
\end{split}
\end{equation*}
for all $n\in\mathcal{T},\,\mathcal{A}\in\mathcal{F}_n,$ which implies 
\begin{equation*}
\begin{split}
      &\Big[b_{u_1}^{*}(n)p_{1}(n)+\sigma_{u_1}^{*}(n)p_{1}(n)\sum_{m=0}^{n-1}d(n,m)\xi_m^H\\&+e(n,n)\sigma_{u_1}^{1,\ast}(n)q_{1}(n)-f_{1,u_1}^{\ast}(n)k_{1}(n)\Big]\\
\cdot&\big(u_1(n)-u_1^{\ast}(n)\big)\geq 0,\quad \text{ a.s. for all $n\in\mathcal{T}$}.
\end{split}
\end{equation*}
In the similar way, we get 
    \begin{equation*}
     \begin{split}
&\Big[b_{u_2}^{*}(n)p_{2}(n)+\sigma_{u_2}^{*}(n)p_{2}(n)\sum_{m=0}^{n-1}d(n,m)\xi_m^H\\&+e(n,n)\sigma_{u_2}^{*}(n)q_{2}(n)-f_{2,u_2}^{\ast}(n)k_{2}(n)\Big]\\
\cdot&\big(u_2(n)-u_2^{\ast}(n)\big)\geq 0, \quad\text{ a.s. for all $n\in\mathcal{T}$}. 
        \end{split}
    \end{equation*}
\begin{rem}
    If the  Nash equilibrium point $\big(u_1^{\ast}(n),u_2^{\ast}(n)\big)$ takes values in the interior of the $\mathbf{U}_1\times \mathbf{U}_2$, then the inequalities in Theorem \ref{19} are equivalent to the following equations
    \begin{equation*}
       \begin{split}
           &b_{u_1}^{*}(n)p_{1}(n)+\sigma_{u_1}^{*}(n)p_{1}(n)\sum_{m=0}^{n-1}d(n,m)\xi_m^H\\&+e(n,n)\sigma_{u_1}^{*}(n)q_{1}(n)-f_{1,u_1}^{\ast}(n)k_{1}(n)=0,a.s.,\\
           & b_{u_2}^{*}(n)p_{2}(n)+\sigma_{u_2}^{*}(n)p_{2}(n)\sum_{m=0}^{n-1}d(n,m)\xi_m^H\\&+e(n,n)\sigma_{u_2}^{*}(n)q_{2}(n)-f_{2,u_2}^{\ast}(n)k_{2}(n)=0,a.s.,
       \end{split} 
    \end{equation*}
    for all $n\in\mathcal{T}$,
    and we denote \begin{equation}\label{eq:H}
\begin{split}
\mathcal{H}&^i\big(n,X(n),u_i(n),u_{-i}^{\ast}(n),p_i(n),q_i(n),k_i(n)\big)=\mathcal{H}^i(n)\\
:=&b^i(n) p_i(n)
+ \sigma^i(n) p_i(n)\sum_{m=0}^{n-1}d(n,m)\xi_m^H
\\
&+ e(n,n)\sigma^i(n) q_i(n)
- f_i(n) k_i(n),
\end{split}
\end{equation}
and 
\begin{equation*}
    \begin{split}
        &\Delta \mathcal{H}^{i}(n)\\=&\mathcal{H}\big(n,X(n),u_i(n),u_{-i}^{\ast}(n),p_i(n),q_i(n),k_i(n)\big) \\&-\mathcal{H}\big(n,X(n),u_i^{\ast}(n),u_{-i}^{\ast}(n),p_i(n),q_i(n),k_i(n)\big) .
    \end{split}
\end{equation*}
\end{rem} 
\begin{thm}[Verification Theorem]\label{thm:verification}
Let \((u_1^{\ast},u_2^{\ast})\in \mathbf{U}_1\times \mathbf{U}_2\) be a Nash equilibrium
point for \textbf{Problem (NZSRG)}, let
\(X^{\ast}=\big\{X^{\ast}(n),\ n\in [0,N]\big\}\) be the corresponding state,
and let
\(\big\{(p_{i}(n),q_{i}(n)),\ n\in [0,N]\big\}\)
be the solution of the adjoint equation \eqref{8}, \(i=1,2\).
Suppose that the following conditions hold:
\begin{enumerate}
\item[(i)] for \(i=1,2\), for all \(n\in [0,N]\) and all
\(u_i\in\mathbf{U}_i\),
\begin{equation*}
    \begin{split}
        &\mathcal{H}^i_{u_i}\big(n,X,u_i^{\ast},u_{-i}^{\ast},p_i,q_i,k_i\big)\big(u_i(n)-u^{\ast}_i(n)\big)\geq 0.
    \end{split}
\end{equation*}
\item[(ii)] $
\mathcal{H}^i\big(n,\cdot,\cdot,u_{-i}^{\ast},p_i,q_i,k_i\big)$ and $g_i(\cdot)$ are convex on \(\mathbf{U}_i\).
\end{enumerate}
Then $(u_1^{\ast},u_2^{\ast})$ is a Nash equilibrium point for
\textbf{Problem (NZSRG)}.
\end{thm}
$Proof.$ For any $u_1(\cdot)\in U_1$, we have
\begin{equation*}\label{16}
\begin{split}
             &\mathbf{E}\big[p_1(N)\big(X^1(N)-X^{*}(N)\big)-p_1(0)\big(X_1(0)-X_1^{\ast}(0)\big)\big]\\+&\mathbf{E}\big[k_1(N)\big(Y_1(N)-Y_1^{\ast}(N)\big)-k_1(1)\big(Y_1(1)-Y_1^{\ast}(1)\big)\big]\\
=&\mathbf{E}\sum_{n=0}^{N-1}\Delta \big[(X^1(n)-X^{*}(n))p_1(n)\big]\\&+\mathbf{E}\sum_{n=1}^{N-1}\big[\Delta\big(Y(n)-Y^{\ast}(n)\big)k_1(n)\big]\\
=&\mathbf{E}\sum_{n=0}^{N-1}\bigg\{\big(X^1(n+1)-X^{*}(n+1)\big)\\&\cdot\Big(-b_x^{*}(n+1)p_1(n+1)\\&-e(n+1,n+1)\sigma_x^{*}(n+1)q_1(n+1)\\
        &+f_{1,x}^{\ast}(n+1)k_1(n+1)\\&-\sigma_x^{*}(n+1)p_1(n+1)\sum_{m=0}^{n-1}d(n,m)\xi_m^H\Big)\\
        &+p_1(n)\big[\big(b^1(n)-b^{*}(n)\big)+\big(\sigma^1(n)-\sigma^{*}(n)\big)\xi_n^H\big]\bigg\}        \end{split}
\end{equation*} 
\begin{equation*}\begin{split}\\&+\mathbf{E}\sum_{n=1}^{N-1}\bigg\{-k_1(n+1)\big[f_1(n+1)-f_1^{\ast}(n+1)\big]\\&+\big(Y_1(n)-Y^{\ast}_1(n)\big)f_{1,y}^{\ast}(n)k_1(n)\\&+\big(Z_1(n)-Z^{\ast}_1(n)\big)f_{1,z}^{\ast}(n)k_1(n)\bigg\}\\ 
=&\mathbf{E}\sum_{n=0}^N\Delta \mathcal{H}^1(n)+\mathbf{E}\sum_{n=1}^N\bigg\{\big(X^1(n)-X^{*}(n)\big)\\&\cdot\Big(-b_x^{*}(n)p_1(n)-e(n,n)\sigma_x^{*}(n+1)q_1(n)\\&+f_{1,x}^{\ast}(n)k_1(n)-\sigma_x^{*}(n)p_1(n)\sum_{m=0}^{n-1}d(n,m)\xi_m^H\Big)\\&+\mathbf{E}\big(Y_1(n)-Y^{\ast}_1(n)\big)\big[f_{1,y}^{\ast}(n)k_1(n)+f_{1,z}^{\ast}(n)k_1(n)\big]\bigg\}        
\\&+\Phi(N),
    \end{split}
\end{equation*}
where
\begin{equation*}
\begin{split}
       & \Phi(N)\\
=&-\mathbf{E}\bigg\{p_1(N)\big[\big(b^1(N)-b^{*}(N)\big)\\&+\big(\sigma^1(N)-\sigma^{*}(N)\big)\xi_n^H\big]\\&+k_1(1)\big(f_1(1)-f_1^{\ast}(1)\big)\\&-\big(Y_1(N)-Y_1^{\ast}(N)\big)f_{1,y}^{\ast}(N)k_1(N)\\&-\big(Z_1(N)-Z_1^{\ast}(N)\big)f_{1,z}^{\ast}(N)k_1(N)\bigg\},
\end{split}
\end{equation*}
and $b^1(N)=b^{*}(N)=\sigma^1(N)=\sigma^{*}(N),$ $f_{1,y}^{\ast}(N)=f_{1,z}^{\ast}(N)=0$.

Then noticing the setting $X(0)=X^{\ast}(0)=k_1(0)$, we get 
\begin{equation*}
    \begin{split}
&J_1\big(u_1(n),u_2^{\ast}(n)\big)-J_1\big(u_1^{\ast}(n),u_2^{\ast}(n)\big)\\=&Y_1(0)-Y_1^{\ast}(0)\\=&\mathbf{E}\big[Y_1(0)+Z_1(0)\eta_n-Y_1^{\ast}(0)-Z_1^{\ast}(0)\eta_n\big]\\=&-\mathbf{E}\big[k_1(1)\big(Y_1(1)-Y_1^{\ast}(1)\big)+k_1(1)\big(f_1(1)-f_1^{\ast}(1)\big)\big]\\
=&\ \mathbf{E}\sum_{n=0}^N\Delta \mathcal{H}^1(n)+\mathbf{E}\sum_{n=0}^N\bigg\{\big(X(n)-X^{\ast}(n)\big)\\&\cdot\Big(-b_x^{*}(n)p_1(n)-e(n,n)\sigma_x^{*}(n+1)q_1(n)\\&+f_{1,x}^{\ast}(n)k_1(n)-\sigma_x^{*}(n)p_1(n)\sum_{m=0}^{n-1}d(n,m)\xi_m^H\Big)\\&+\mathbf{E}\big(Y_1(n)-Y^{\ast}_1(n)\big)\big[f_{1,y}^{\ast}(n)k_1(n)+f_{1,z}^{\ast}(n)k_1(n)\big]\bigg\}        
\\\geq &\mathbf{E}\sum_{n=0}^N \mathcal{H}^1_{u_1}(n,X^{*},u_i^{\ast},u_{-i}^{\ast},p_i,q_i,k_i)\big(u_1(n)-u_1^{\ast}(n)\big)\\\geq& 0.
\end{split}
\end{equation*}
Therefore, we prove the theorem.

\subsection{Maximum principle for discrete-time zero-sum stochastic game}
In this subsection,
we consider the system with (\ref{3}), (\ref{y_0}) and (\ref{26}).
 We still assume that $b,\sigma,f_i,g_i$ satisfy the assumptions (H1)--(H4). 

\textbf{Problem (NZSRG1)} Denote $(f_1,g_1)=(f,g)$ and $(f_2,g_2)=(-f,-g)$. Then the saddle point $(u_1^{\ast},u_2^{\ast})$ of \textbf{Problem (ZSRG)} satisfies 
\begin{equation*}
    \begin{split}
&J_1(u_1^{\ast},u_2^{\ast})=\inf_{u_1\in U_1}J_1(u_1,u_2^{\ast}),\\&J_2(u_1^{\ast},u_2^{\ast})=\inf_{u_2\in U_2}J_2(u_1^{\ast},u_2),
    \end{split}
\end{equation*}
and 
\begin{equation*}
J_i(u_1,u_2)=Y_i(0),\quad\quad i=1,2,
\end{equation*}
where $\big(Y_i(n),Z_i(n)\big)$ is the solution to the BS$\Delta$E (\ref{eq:adjoint}).

Note that the problem defined by system (\ref{3}) with cost functionals  $J_1$  and  $J_2$  constitutes a nonzero-sum game of the type introduced in \textbf{Problem (NZSRG)}. Consequently, \textbf{Problem (ZSRG)} can be transformed into \textbf{Problem (NZSRG1)}. Furthermore, a saddle point $(u_1^*, u_2^*)$ for \textbf{Problem (ZSRG)} corresponds precisely to a equilibrium point for \textbf{Problem (NZSRG1)}. The necessary condition for the saddle point is then derived based on these definitions.
\begin{thm}\label{47}
    Suppose $(u_1^{\ast},u_2^{\ast})$ is a saddle point of \textbf{Problem (ZSRG)} and $X^{\ast}(n)$ is the state variable. Then for each $n\in\mathcal{T}$ , it holds that 
    \begin{equation}\label{20}
        \begin{split}
            \begin{cases}       
            \mathcal{H}_{u_1}^{1,*}(n)\cdot\big(u_1(n)-u_1^{\ast}(n)\big) \ge 0,\\ \mathcal{H}_{u_2}^{2,\ast}(n)\cdot\big(u_2(n)-u_2^{\ast}(n)\big) \le 0.
\end{cases}
        \end{split}
    \end{equation}
   where $\big(p(n),q(n),k(n)\big)$ is the solution of the following adjoint equation 
    \begin{equation}\label{k}
        \begin{split}
            \begin{cases}
                k(n+1)=k(n)+f_y^{\ast}(n)k(n)+f_z^{\ast}(n)k(n)\eta_n,\\p(n)+q(n)\eta_n=p(n+1)+H_x^{\ast}(n,X,u_i,u_{-i},p,q,k)\\
                k(1)=-1,\quad k(0)=0,\\
                p(N)=-g_x\big(X^{\ast}(N)\big)k(N).
            \end{cases}
        \end{split}
    \end{equation}
\end{thm}

    $Proof.$ By Theorem \ref{19}, we have for each $n\in\mathcal{T}$ , 
    \begin{equation}\label{18}
        \begin{split}
            \begin{cases} 
            \mathcal{H}_{u_1}^{1,\ast}(n)\cdot\big(u_1(n)-u_1^{\ast}(n)\big) \ge 0,\\ \mathcal{H}_{u_2}^{2,\ast}(n)\cdot\big(u_2(n)-u_2^{\ast}(n)\big) \ge 0.
            \end{cases}
        \end{split}
    \end{equation}
By the definition of $\mathcal{H}^i$, it's easy to check that $\big(p_{i}(n),q_{i}(n),k_{i}(n)\big)$ satisfies
 \begin{equation*}
        \begin{split}
            \begin{cases}
                k_i(n+1)=k_i(n)+f_{i,y}^{\ast}(n)k_i(n)+f_{i,z}^{\ast}k_i(n)\eta_n,\\p_i(n)+q_i(n)\eta_n=p_i(n+1)+b_x^{*}(n+1)p_i(n+1)\\\quad\quad\quad\quad\quad\quad\quad+e(n+1,n+1)\sigma_x^{*}(n+1)q_i(n+1)\\\quad\quad\quad\quad\quad\quad\quad +\sigma_x^{*}(n+1)p_i(n+1)\mathbf{E}^{\mathcal{F}_{n+1}}\big[\xi_{n+1}^H\big]\\\quad\quad\quad\quad\quad\quad\quad-f_{i,x}^{\ast}(n+1)k_i(n+1),\\
                k_i(1)=1,\quad k_i(0)=0,\\ p_i(N)=-g_{i,x}\big(X^{\ast}(N)\big)k_i(N),\quad\quad i=1,2.
            \end{cases}
        \end{split}    \end{equation*}    Comparing these two equations with the case of $i=1$ and the case of $i=2$, we get 
        $\big(p_{1}(n),q_1(n),k_1(n)\big)=\big(-p_2(n),-q_2(n),k_2(n)\big)$ for each $n\in\mathcal{T}$. 
        Define $\big(p(n),q(n),k(n)\big)=\big(p_1(n),q_1(n),k_1(n)\big)$ for each $n\in\mathcal{T}$. Finally, $\mathcal{H}^2=-\mathcal{H}^1$, it's easy to check that (\ref{18}) is just (\ref{20}).
\par Next result is the sufficient condition of the saddle point.

\begin{thm}[Verification Theorem]
\label{thm:zsg-sufficient-recursive}
Let \((u_1^*,u_2^*)\in U_1\times U_2\) be a pair of admissible controls for \textbf{Problem (ZSRG)}, and let \(X^*=\{X^*(n),n\in\mathcal{T}\}\) be the corresponding state process. Let $\big(p(n),q(n),k(n)\big)$ satisfies (\ref{k}).
Assume that for each \(n\in \mathcal{T}\),
\(\mathcal{H}\big(n,\cdot,\cdot,u_2^{\ast},p,q,k\big)\) is convex,
\(\mathcal{H}\big(n,\cdot,u_1^{\ast},\cdot,p,q,k\big)\) is concave, and
\(g(x)=Wx+\gamma\), where \(W\) is a bounded \(\mathcal{F}_N\)-measurable random variable and \(\gamma\) is an \(\mathcal{F}_N\)-measurable integrable random variable.
Then \((u_1^*,u_2^*)\) is a saddle point of Problem (ZSG) if (\ref{20}).
\end{thm}
$Proof.$ 
The proof is separated into three parts.
\par
 Part 1. Since $(g_1,f_1) = (g,f)$ and $\big(p_1(n),q_1(n),k_1(n)\big)$ $= \big(p(n),q(n),k(n)\big)$, $\mathcal{H}\big(n,\cdot,\cdot,u_{2}^{\ast}(n),p(n),q(n)\big)$ and $g(\cdot)$ are convex. We have $J_1(u_1^{\ast},u_2^{\ast}) = \inf _{ u_1\in\mathbf{U}_1} J_1(u_1,u_2^{\ast})$. That is, $J(u_1^{\ast}, u_2^{\ast}) = \inf_{u_1\in\mathbf{U}_1} J(u_1, u_2^{\ast}).$

Part 2. Since $(g_2,f_2) = (-g,-f)$ and $\big(-p_2(n),-q_2(n)$ $,k_2(n)\big)= \big(p(n),q(n),k_2(n)\big)$, $\mathcal{H}\big(n,\cdot,u_{1}^{\ast}(n),\cdot,p(n),q(n)\big)$ and $g(\cdot)$ are convex.  
$J_2(u_1^{\ast},u_2^{\ast}) = \inf _{ u_2\in\mathbf{U}_2} J_2(u_1,u_2^{\ast})$. That is, $J(u_1^{\ast}, u_2^{\ast}) = \sup_{u_2\in\mathbf{U}_2} J(u_1, u_2^{\ast}).$
\par Part 3. Since $x\rightarrow g(x)$ is assumed to be linear, the conclusion is a combination of Parts 1 and 2.

\section{Applications}

\subsection{Discrete-time LQ nonzero-sum stochastic game}
Let us study discrete-time LQ nonzero-sum stochastic games. Take $U_{i}(n)=\mathbf{R}^m$ for all $i=1,2$ and $n\in\mathcal{T}$. Assume $d=1.$

\textbf{Problem (LQNZSRG)} The problem is to find $(u_1^{\ast},u_2^{\ast})\in\mathbf{U}_1\times\mathbf{U}_2$ such that (\ref{31}) holds, the cost functional is given as 
\begin{equation*}
J_i(u_1,u_2)=Y_i(0),
\end{equation*}
where
\begin{equation*}
\begin{split}
\begin{cases}
    Y_{i}(n)+Z_{i}(n)\eta_n=Y_{i}(n+1)+G_{i}(n+1)X(n+1)\\\quad\quad\quad\quad\quad\quad\quad\quad+(H_{i}(n+1)-1)Y_i(n+1)\\\quad\quad\quad\quad\quad\quad\quad\quad+\frac{1}{2}L_{1}(n+1)u_{1}^2(n+1)\\\quad\quad\quad\quad\quad\quad\quad\quad+\frac{1}{2}L_{2}(n+1)u_{2}^2(n+1),\quad n\in\mathcal{T},\\
    Y_{i}(N)=R_{i}(N)X(N), 
\end{cases}
    \end{split}
\end{equation*}
and 
\begin{equation}\label{35}
    \begin{split}
        \begin{cases}
X(n+1)=X(n)+\big(A(n)X(n)+B_{1}(n)u_{1}(n)\\\quad\quad\quad\quad\quad+B_{2}(n)u_{2}(n)+C(n)\big)+\big(D(n)X(n)\\
\quad\quad\quad\quad\quad +E_{1}(n)u_{1}(n)+E_{2}(n)u_{2}(n)\\\quad\quad\quad\quad\quad+F(n)\big)\xi_n^H, \quad n\in\mathcal{T}\\
X(0)=x_0 .
        \end{cases}
    \end{split}
\end{equation}
The coefficients are deterministic and bounded matrices of appropriate dimensions. Moreover, $G_{i}(n),L_i>0,R_{i}(n)\geq 0$ for $i=1,2.$ 
According to the section 3, the adjoint equations are as follows,
\begin{equation*}
    \begin{split}
        \begin{cases}
            k_i(n)=-\Pi _{j=1}^{n-1}H_i(j),\quad n=1,\cdots,N,\\
            p_{i}(n)+q_{i}(n)\eta_n=p_{i}(n+1)+A(n+1)p_{i}(n+1)\\\quad\quad\quad\quad\quad\quad+e(n+1,n+1)D(n+1)q_{i}(n+1)\\
            \quad\quad\quad\quad\quad\quad+D(n+1)p_{i}(n+1)\mathbf{E}^{\mathcal{F}_{n+1}}\big[\xi_{n+1}^H\big]\\\quad\quad\quad\quad\quad\quad-G_{i}(n+1)k_{i}(n+1),\quad n\in\mathcal{T},\\k_i(1)=-1,\quad k_i(0)=0\\
            p_{i}(N)=-R_{i}(N)k_i(N).
        \end{cases}
    \end{split}
\end{equation*}
\par The necessary conditions turn out to be 
\begin{equation*}
    \begin{split}
        \begin{cases}
            u_{1}^{\ast}(n)=-\big[L_{1}(n)\Pi_{j=1}^{n-1}(j)\big]^{-1}\big(B_{1}(n)p_{1}(n)\\\quad\quad\quad\quad+e(n,n)E_{1}(n)q_{1}(n)\\\quad\quad\quad\quad+E_{1}(n)p(n)\sum_{m=0}^{n-1}d(n,m)\xi_k^H\big),\\
              u_{2}^{\ast}(n)=-\big[L_{2}(n)\Pi_{j=1}^{n-1}(j)\big]^{-1}\big(B_{2}(n)p_{2}(n)\\\quad\quad\quad\quad+e(n,n)E_{2}(n)q_{2}(n)\\\quad\quad\quad\quad+E_{2}(n)p(n)\sum_{m=0}^{n-1}d(n,m)\xi_k^H\big),
        \end{cases}
    \end{split}
\end{equation*}
where $n\in\mathcal{T}$.
In the remaining part of this example, we consider the simple case when $\big\{D(n)\big\},\big\{E_{1}(n)\big\}$ and $\big\{E_{2}(n)\big\}$ vanish. Then a state-depending form of the equilibrium point $(u_1^{\ast},u_2^{\ast})$ is given. The noises and admissible controls are allowed to be only square-integrable. In this case, it holds that
\begin{equation}\label{32}
     \begin{split}
        \begin{cases}
         u_{1}^{\ast}(n)=-\big[L_{1}(n)\Pi_{j=1}^{n-1}(j)\big]^{-1}\big(B_{1}(n)p_{1}(n)\big),\\u_{2}^{\ast}(n)=-\big[L_{2}(n)\Pi_{j=1}^{n-1}(j)\big]^{-1}\big(B_{2}(n)p_{2}(n)\big).\end{cases}
    \end{split}
\end{equation}
Iterating backward yields, we can get the explicit solution of $p_i(n)$ as follows,
\begin{equation*}
\begin{split}
       p_i(n)
={}&R_i(N)
\bigg(\prod_{k=n+1}^{N}(1+A(k))\bigg)
\bigg(\prod_{j=1}^{N-1}H_i(j)\bigg)\\
+\sum_{k=n+1}^{N}&
\bigg[
G_i(k)
\bigg(\prod_{\ell=n+1}^{k-1}(1+A(\ell))\bigg)
\bigg(\prod_{j=1}^{k-2}H_i(j)\bigg)
\bigg].
\end{split}
\end{equation*} 
According to the existence and uniqueness results of the adjoint equation, $q_i(n) = 0$.

Then by equation (\ref{32}), it holds that 
\begin{equation}\label{34}
     \begin{split}
            &u_{i}^{\ast}(n)\\=&-\Big[L_{i}^{-1}(n)\Pi_{j=1}^{n-1}(j)\Big]\\&\cdot\Bigg\{R_i(N)
\bigg(\prod_{k=n+1}^{N}(1+A(k))\bigg)
\bigg(\prod_{j=1}^{N-1}H_i(j)\bigg)\\
&+\sum_{k=n+1}^{N}
\bigg[
G_i(k)
\bigg(\prod_{\ell=n+1}^{k-1}(1+A(\ell))\bigg)
\bigg(\prod_{j=1}^{k-2}H_i(j)\bigg)
\bigg]\Bigg\}.
    \end{split}
\end{equation}
\subsection{Discrete-time LQ zero-sum stochastic game}
For a linear-quadratic zero-sum stochastic differential game with recursive cost functional (\textbf{Problem (LQZSRG)}), the Nash equilibrium is derived using the stochastic maximum principle. Below is a concise derivation and the explicit form of the Nash equilibrium.

Consider the following state equation,
\begin{equation*}
\begin{split}
    \begin{cases}
         X(n+1) =X(n)+ A(n)X(n) + B_1(n)u_1(n) \\\quad\quad\quad\quad\quad+ B_2(n)u_2(n) +\big( C(n)X(n)+S(n)\big) \xi_n^H, \\X(0) = x_0,\quad n\in [0,N]  
    \end{cases}
\end{split}
\end{equation*}
and the recursive BS$\Delta$E,
\begin{equation*}
\begin{split}
    \begin{cases}
        Y(n)+Z(n)\eta(n)= Y(n+1)+\big[ Q(n+1)Y(n+1) \\\quad\quad\quad\quad\quad\quad\quad\quad- R_1(n+1)u_1^2(n+1)\\\quad\quad\quad\quad\quad\quad\quad\quad+R_2(n+1)u_2^2(n+1) \big] ,\\ Y(N) =  G X(N),\quad n\in\mathcal{T}.
    \end{cases}
\end{split}
\end{equation*}
The cost functional is 
$J(u_1,u_2)=Y(0)$,
Player 1 minimizes $J$, while Player 2 maximizes $J$. We seek a saddle point $(u_1^{\ast}, u_2^{\ast})$ satisfying,
\begin{equation*}
    J(u_1^*, u_2) \leq J(u_1^*, u_2^*) \leq J(u_1, u_2^*), \forall (u_1,u_2)\in\mathbf{U}_1\times\mathbf{U}_2.
\end{equation*}
According to the section 3.2, we can transform the (\textbf{Problem (ZSRG)}) into (\textbf{Problem (NZSRG)}). 
Here $X(n) \in \mathbf{R}$ is the state (one-dimensional for simplicity, but extendable to higher dimensions). $u_1(n) $ is the control of Player 1 (the minimizer). $u_2(n) $ is the control of Player 2 (the maximizer). All coefficients $A, B_i,C,M,Q_i, R_i, G_i$ are deterministic and continuous, with $R_i(n) > 0,i=1,2$.

If the $(p_i, q_i, k_i)$ satisfies the two conditions of Theorem \ref{47}, then we have
\begin{equation*}
\begin{split}
    \begin{cases}
         p_i(n)+q_i(n)\eta_n=p_i(n+1) +A(n) p_i(n+1) \\\quad\quad\quad\quad\quad\quad\quad\quad+e(n+1,n+1) C(n) q_i(n+1)\\\quad\quad\quad\quad\quad\quad\quad\quad+C(n+1)p_i(n+1)\mathbf{E}^{\mathcal{F}_{n+1}}\big[\xi_{n+1}^H\big], 
         \\ 
k_i(n+1) = k_i(n)+Q_ik_i(n),\quad n\in\mathcal{T}, \\ p_i(N) = -k_i(N) G_i X(N), \\k_i(1) = -1,\,k_i(0)=0,\quad i=1,2.
    \end{cases}
\end{split}
\end{equation*}
From the Theorem \ref{47}, the Nash equilibrium point $\big(u_1^{\ast}(n),u_2^{\ast}(n)\big)$ is as follows, 
\begin{equation*}
    \begin{split}
    u_1^{\ast}(n) = -\frac{p_1(n) B_1(n)}{2 R_1(n)k_1(n)}, u_2^{\ast}(n) = \frac{p_2(n) B_2(n)}{2 R_2(n)k_2(n)}.
    \end{split}
\end{equation*}

We assume a linear relationship between the adjoint variable $p$ and the state $X$,
$p_i(n) = P_i(n) X(n),$
where $P_i(n)$ is a deterministic function. Substituting into the expressions for $u_1^{\ast}(n)$ and $u_2^{\ast}(n)$ yields the state feedback form,
\begin{equation*}
\begin{split}
      u_1^{\ast}(n) =& -\frac{P_1(n) B_1(n)}{2 R_1(n)k_1(n)}X(n),\\u_2^{\ast}(n) =& \frac{P_2(n) B_2(n)}{2 R_2(n)k_2(n)}X(n), 
\end{split}
\end{equation*}
where $P_i(n)$ is the solution to the Riccati equation as follows.
To determine $P_i(n)$, substitute $p_i(n) = P_i(n) X(n)$ and $u_i^{\ast}(n)$ into the state equation. In the remaining part of this example, we consider the simple case when $C(n)$ and $S(n)$ vanish. We obtain the Riccati equation for $P_i(n)$,
\begin{equation*}
    \begin{split}
        \begin{cases}
            P_i(n)=P_i(n+1)\big(1+A(n)\big)\\
            \quad\quad\quad\quad\cdot\big[1+A(n)-\frac{P_1(n)B_1^2(n)}{2R_1(n)k_1(n)}+\frac{P_2(n)B_2^2(n)}{2R_2(n)k_2(n)}\big], \\ P_i(N) = -k_i(N)G_i,\quad i=1,2.
        \end{cases}
    \end{split}
\end{equation*}

Next, we simplify the coupled Riccati Equations.
 
\textbf{Step 1}, both equations share the same multiplicative factor on the right-hand side. 
\begin{equation*}
\begin{split}
    P_i(n) =& P_i(n+1) \cdot \Phi(n) \\&\cdot\big[ \Phi(n) + \Psi_1(n) P_1(n) + \Psi_2(n) P_2(n) \big],
\end{split}
\end{equation*}
where
\begin{equation*}
\begin{split}
     \Phi(n) =& 1 + A(n), \\\Psi_1(n) =& -\frac{B_1^2(n)}{2R_1(n)k_1(n)}, \Psi_2(n) = \frac{B_2^2(n)}{2R_2(n)k_2(n)}.
\end{split}
\end{equation*}
Since the bracketed term is identical for both $i$, the ratio $P_1(n)/P_2(n)$ is constant for all $n$. Let $ \alpha = \frac{P_1(n)}{P_2(n)}$.

From the terminal conditions,
\begin{equation*}
    \alpha = \frac{P_1(N)}{P_2(N)} = \frac{-k_1(N)G_1}{-k_2(N)G_2} = \frac{k_1(N)G_1}{k_2(N)G_2}.
\end{equation*}

\textbf{Step 2}, 
substitute $P_1(n) = \alpha P_2(n)$ into the equation as follows,
\begin{equation*}
\begin{split}
     P_2(n) =& P_2(n+1) \Phi(n) 
     \\&\cdot\left[ \Phi(n) + \Psi_1(n) \alpha P_2(n) + \Psi_2(n) P_2(n) \right].
\end{split}
\end{equation*}
Define
\begin{equation*}
\begin{split}
    \Psi(n) =& \Psi_1(n) \alpha + \Psi_2(n) \\=& -\frac{B_1^2(n)}{2R_1(n)k_1(n)} \alpha + \frac{B_2^2(n)}{2R_2(n)k_2(n)}.
\end{split}
\end{equation*}

\textbf{Step 3}, to solve for $P_2(n)$, invert by setting $\Gamma(n) = 1 / P_2(n)$ (assuming $P_2(n) \neq 0$). Substituting,
\begin{equation*}
\begin{split}
     \frac{1}{\Gamma(n)} =& \frac{P_2(n+1) \Phi^2(n)}{1 - P_2(n+1) \Phi(n) \Psi(n)}\\ = &\frac{\Phi^2(n) / \Gamma(n+1)}{1 - \Phi(n) \Psi(n) / \Gamma(n+1)} \\=& \frac{\Phi^2(n)}{\Gamma(n+1) - \Phi(n) \Psi(n)},
\end{split}
\end{equation*}
which yields
\begin{equation*}
    \Gamma(n) = \frac{\Gamma(n+1) - \Phi(n) \Psi(n)}{\Phi^2(n)}.
\end{equation*}
This is a linear backward recurrence for $\Gamma(n)$.

\textbf{Step 4}, the recurrence
$\Gamma(n) = \frac{\Gamma(n+1)}{\Phi^2(n)} - \frac{\Psi(n)}{\Phi(n)}$, with terminal condition,
\begin{equation*}
    \Gamma(N) = \frac{1}{P_2(N)} = -\frac{1}{k_2(N) G_2}.
\end{equation*}
The explicit solution is
\begin{equation*}
   \Gamma(n) = \Phi(n, N) \, \Gamma(N) - \sum_{k=n}^{N-1} \Phi(n, k+1) \frac{\Psi(k)}{\Phi(k)}, 
\end{equation*}
where $\Phi(n, m)$ is the transition function for the homogeneous part,
\begin{equation*}
    \Phi(n, m) = \prod_{j=n}^{m-1} \frac{1}{\Phi^2(j)}, \quad \text{for } n < m, \quad \text{and} \quad \Phi(n, n) = 1.
\end{equation*}
Then
\begin{equation*}
   P_2(n) = \frac{1}{\Gamma(n)}, \quad P_1(n) = \alpha P_2(n). 
\end{equation*}

\subsection{Portfolio Selection Problem}
Similar to the works of Ji and Zhang \cite{Yuanyuan}, Song and Liu \cite{Teng} and Browne \cite{Siddifferentialportfoliogames}, this paper examines a portfolio selection problem involving a recursive utility functional within a financial market. Drawing on the framework introduced in \cite{Espinosa2015OptimalInvestmen}, we incorporate a relative performance component into the analysis. This component shifts the focus from a conventional optimization problem to a game-theoretic setting, where each investor's utility depends not only on their own wealth but also on that of other investors, thereby introducing strategic interactions. We consider a financial market model consisting of two investment opportunities.
\par One is a risk-free asset with the interest rate $r$, such as a bank deposit. The total amount invested in the risk-free asset is $S_0$,
\begin{equation*}
    \begin{split}
        \begin{cases}
            S_0(n+1)-S_0(n)=r(n)S_0(n),\\
            S_0(0)=s_0>0,
        \end{cases}
    \end{split}
\end{equation*}
where $r(n)$ is a bounded deterministic function.
\par The other is risky asset with the return rate $\mu$ and the volatility $\sigma$, such as stocks.  and the total amount invested in the risky assets is $S_1$. Assume that $S_0$ and $S_1$ satisfy the following equations,
\begin{equation*}
    \begin{split}
    \begin{cases}
       S_{1}(n+1)-S_{1}(n)=\mu(n)S_{1}(n)+\sigma(n)S_{1}(n)\xi_n^H,\\
       S_{1}(0)=s_1>0,i=1,2,
    \end{cases}
    \end{split}
\end{equation*}
Where $\mu(n),\sigma(n)\neq 0$ are bounded deterministic functions, $\mu(n)>r(n)>0,$ and $\xi_n^H$ is a Gaussian process with long-range dependence.
The total corresponding wealth process $X(n)$ is formulated as 
\begin{equation*}
    X(n)=\theta_0(n)S_0(n)+\theta_1(n)S_1(n),\quad n\in\mathcal{T},
\end{equation*}
where $\theta_0(n)$ and $\theta_1(n)$ represent the investment amounts in risk-free assets and risky assets, respectively.
For a self-financing portfolio, the difference in total wealth $X$ between two adjacent discrete times is purely owing to the change in the prices of the stocks. Hence , using the self-financing condition, see Song and Liu \cite{Teng}. By adding the self-financing conditions of the two investors, the self-financing condition of the total wealth is obtained,
\begin{equation*}
\begin{split}
     &\theta_{0}(n+1)S_0(n+1)+\theta_{1}(n+1)S_1(n+1)\\=&\theta_{0}(n)S_0(n+1)+\theta_{1}(n)S_1(n+1),
\end{split}
\end{equation*}
and let 
\begin{equation*}
    \theta_0(n)=(\theta_{0,1}(n),\theta_{0,2}(n)),\theta_1(n)=(\theta_{1,1}(n),\theta_{1,2}(n)).
\end{equation*}

Denote $u_n$ the amount invested in the risky security, so we can get the risk asset holdings and risk-free asset holdings
\begin{equation*}
    \begin{split}
        \theta_{1,i}(n)=\dfrac{u_i(n)}{S_{1,i}(n)},\quad \theta_{0,i}(n)=\dfrac{X_i(n)-u_i(n)}{S_{0,i}(n)}.
    \end{split}
\end{equation*}
According to the above equations, we obtain the wealth dynamics
\begin{equation*}
    \begin{split}
    \begin{cases}
X(n+1)=X(n)+r(n)X(n)+\big(\mu(n)-r(n)\big)\\\quad\quad\quad\quad\cdot\big(u_1(n)+u_2(n)\big)+\sigma(n)\big(u_1(n)+u_2(n)\big)\xi_n^H,\\
       X(0)=x_{0}>0,\quad i = 1,2.
    \end{cases}
    \end{split}
\end{equation*}
We denote by  $\mathbf{U}_{i}$  the set of admissible portfolios taking values in  $\mathbf{R}$. Consider a small investor with initial wealth $x_0$, who selects a portfolio strategy $u_i(n)$ at each time $ n$. The investor aims to choose an optimal strategy $u_i^{\ast}(\cdot) \in \mathbf{U}_{i}$ that maximizes the expected utility functional with reward function  
\begin{equation*}
    \begin{split}
    \begin{cases}
        f_i(n) = -C(n)X(n)+D_i(n)u_i^{\beta_i}(n)+(\alpha_i-1) Y_i(n),\\
  g_i(X(N))=-R_iX(N),\quad\quad i=1,2, 
    \end{cases}
    \end{split}
\end{equation*}
where $x(n)$ represents the total market wealth or a common state process. 
 $u_i(n)$ is the amount invested in the risky asset by investor  $i$  at time  $n$. Parameters $\beta_i > 1$ are related to risk preference. $\alpha_i\in(0,1)$  is the discount factor. $ R_i > 0$ represents weight for terminal wealth. The negative sign in the terminal function  $g_i(X(N)) = -R_i x(N)$  indicates that the investor aims to maximize terminal wealth.

Denote $u_i^{\ast}(\cdot)$ be an optimal portfolio strategy and $X^{\ast}(\cdot),Y^{\ast}(\cdot)$ be the corresponding wealth process, recursive utility process, respectively, with corresponding solution $(p^{\ast}(\cdot),q^{\ast}(\cdot))$ of the adjoint equation. Then we can use the Theorem \ref{19} to get the adjoint equations as follows, 
\begin{equation*}
    \begin{split}
        \begin{cases}
            k_{i}(n)= -\alpha_i^{n-1}, \\
            p_i(n)+q_i(n)\eta_n=p_i(n+1)\\\quad\quad\quad\quad\quad\quad\quad\quad+\big[r(n)p_i(n+1)+C(n)k_{i}(n+1)\big],\\
            p_i(N)=R_ik_i(N),
            k_1=-1,\quad \quad n\in\mathcal{T}.
        \end{cases}
    \end{split}
\end{equation*}

Therefore, we have the following explicit expressions for $p_i(n)$ and $q_i(n)$,
    \begin{equation*}
       \begin{split}
       \begin{cases}
              p_i(n) = \bigg(\prod_{j=n+1}^{N} \big(1 + r(j)\big)\bigg)p_N\\\quad\quad\quad-\sum_{m=n+1}^N\bigg(\prod_{j=n+1}^{m-1} \big(1 + r(j)\big)\bigg)C(m)\alpha_i^{m-1}, \\
            q_i(n) = 0,\quad i=1,2, \quad n\in\mathcal{T}.
       \end{cases}
       \end{split}
   \end{equation*}
    If the Nash equilibrium point $\big(u^{\ast}_1(n),u_2^{\ast}(n)\big)$ takes values in the interior of the $\mathbf{U}_1\times\mathbf{U}_2$, the Nash equilibrium point is given by   
    \begin{equation*}
        \begin{split}
    u^{\ast}_{i}(n)=\bigg[\big(\mu(n)-r(n)\big)p_i(n)+\sigma(n)p_i(n)\sum_{m=0}^{n-1}d(n,m)\xi_m^H\\+e(n,n)\sigma(n)q_i(n)\bigg]^{\frac{1}{\beta_i-1}}\bigg[\alpha_i^{n-1}D_i(n)\beta_i\bigg]^{-\frac{1}{\beta_i-1}},n\in\mathcal{T}.  
        \end{split}
    \end{equation*}
    
  In the rest part of this example, we take the simple case as follows. Here we assume that $r,b,\sigma,C,D_i,\alpha_i,R_i$ are deterministic constant or process, so $p_i(n)$ is also deterministic. The explicit expressions for $p_i(n)$ and $q_i(n)$ are as follows,
     \begin{equation*}
       \begin{split}
       \begin{cases}
              p_i(n) = \bigg( \big(1 + r\big)^{N-n}\bigg)p_N\\\quad\quad\quad-\sum_{j=n+1}^N\bigg(\big(1 + r\big)^{j-n-1}\bigg)C\alpha_i^{j-1}, \\
            q_i(n) = 0,\quad i=1,2, \quad  n\in\mathcal{T}.
       \end{cases}
       \end{split}
   \end{equation*}
   Then, the  Nash equilibrium point satisfies the two equations as follows,
    \begin{equation*}
        \begin{split}
    u^{\ast}_{i}(n)
    =&\bigg\{\big[(\mu-r)+\sigma \sum_{m=0}^{n-1}d(n,m)\xi_m^H\big]\\&\cdot\big[R_1\big((1+r)\alpha_1\big)^{N-n}\\&+\sum_{j=n+1}^N\big((1+r)^{j-n-1}\big)C\alpha_1^{j-n}\big]\bigg\}^{\frac{1}{\beta_i-1}}\\&\cdot\bigg\{{D_i\beta_i}\bigg\}^{-\frac{1}{\beta_i-1}},\quad\quad\quad\quad i=1,2,\quad n\in\mathcal{T}.  
        \end{split}
    \end{equation*}
 \begin{rem}
      The formation of a game occurs because investor $i$'s investment decision $u_i(k)$ directly affects the total wealth $X(k+1),$ thereby influencing the wealth of all investors.
  \end{rem}
  
Finally, we show some numerical results for the Nash equilibrium $\{u_i^{\ast}(n), n\in\mathcal{T}\}$.

In Figure 1, we take the following parameter configuration:  $H = 0.6$,  $\mu = 0.08$, $\sigma = 0.8$ , $r = 0.06$,  $C = 0.12$,  $D_1 = 1$, $D_2 = 1.15$, $\alpha_1 = 0.85$, $\alpha_2 = 0.8$, $R_1 = R_2 = 20$, $N=30$. Figure 1 plots the resulting equilibrium strategy profiles $u_i^*(n)$.
Two primary comparisons are illustrated. Firstly, a cross‑investor comparison (Player 1 versus Player 2) under risk‑aversion coefficients $\beta$ reveals their strategic interaction, the difference in $\beta$ generates distinct investment patterns that reflect the competitive nature of the underlying differential game. Secondly, a within‑player comparison for Player 1 across three panels demonstrates the sensitivity of the optimal strategy to its own risk‑aversion parameter. As $\beta_1$ increases from 2 (Panel a) to 3 (Panel b) and further to 5 (Panel c), $u_1^*(n)$ shifts toward a more aggressive allocation to the risky asset, indicating a increase in risk‑taking.

\begin{figure}[htbp]
\centering
\begin{minipage}{\linewidth} 
    \centering
    \includegraphics[width=0.7\linewidth]{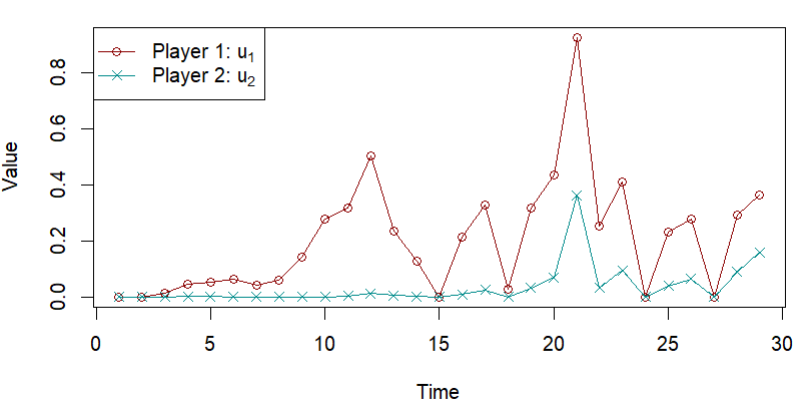} 
    \subfigure{$(a)\quad\beta_1=2,\beta_2=1.5$} 
    \label{fig:group1_a}
\end{minipage}

\vspace{0.2cm} 

\begin{minipage}{\linewidth}
    \centering
    \includegraphics[width=0.7\linewidth]{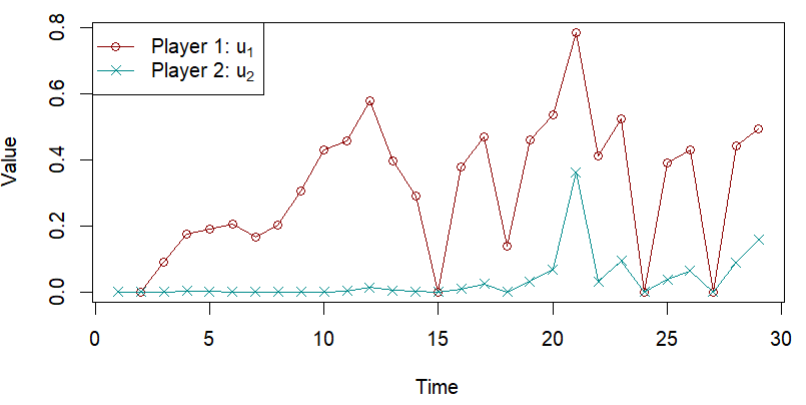}
    \subfigure{$(b)\quad\beta_1=3,\beta_2=1.5$} 
    \label{fig:group1_b}
\end{minipage}

\vspace{0.2cm}

\begin{minipage}{\linewidth}
    \centering
    \includegraphics[width=0.7\linewidth]{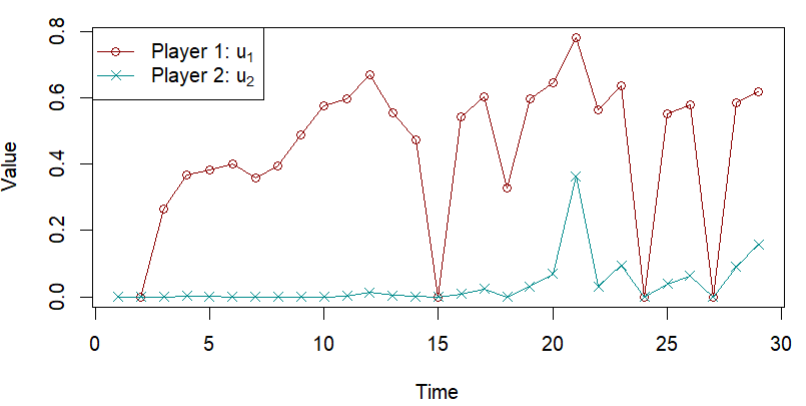}
    \subfigure{$(c)\quad\beta_1=5,\beta_2=1.5$} 
\label{ }
\end{minipage}

\caption{Trajectories $u_1(n)$ vs $u_2(n)$ under different risk preference $\beta$}
\label{fig:group1}
\end{figure}

In Figure 2, with the parameter set $\mu = 0.2, r = 0.07, \sigma = 0.5, C = 0.12, D_1 = 1.05, D_2 = 1, \alpha_1 = 0.85, \alpha_2 = 0.8, R_1 = 20, R_2 = 20, \beta_1 = 2, \beta_2 = 2,$ and $N=30$, the resulting strategy paths are displayed. If $H$ equals $1/2$, the process follows standard Brownian motion and appears relatively smooth.
These comparative plots why contemporary research typically focus on the case $H \neq 1/2$. If $H =1/2$, its smooth trajectory fails to capture the facts commonly observed in real financial markets, such as volatility clustering and long-range dependence.

\begin{figure}[htbp]
\centering
\begin{minipage}{\linewidth} 
    \centering
    \includegraphics[width=0.75\linewidth]{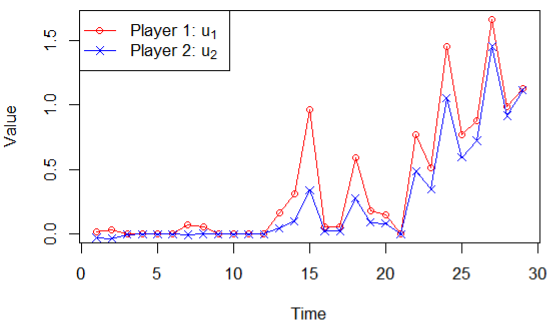} 
    \subfigure{$(a)\quad H=0.3$} 
    \label{fig:group1_a}
\end{minipage}

\vspace{0.2cm} 

\begin{minipage}{\linewidth}
    \centering
    \includegraphics[width=0.75\linewidth]{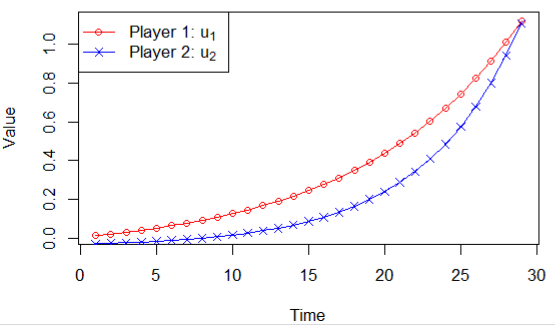}
    \subfigure{$(b)\quad H=0.5$} 
    \label{fig:group1_b}
\end{minipage}

\vspace{0.2cm}

\begin{minipage}{\linewidth}
    \centering
    \includegraphics[width=0.75\linewidth]{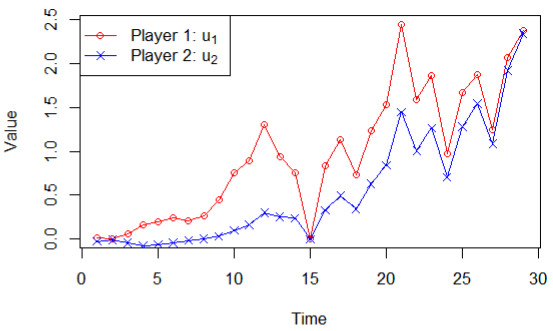}
    \subfigure{ $(c)\quad H=0.8$} 
    \label{fig:group1_c}
\end{minipage}

\caption{Trajectories $u_1(n)$ and $u_2(n)$ with different Hurst exponent $H$}
\label{fig:group1}
\end{figure}

In Figure 3, we take $H = 0.6, \sigma = 0.8, C = 0.12, D_1 = 1.05, D_2 = 1, \alpha_1 = 0.85, \alpha_2 = 0.8, R_1 = 20, R_2 = 20, \beta_1 = 2, \beta_2 = 2, N=30$. Figure 3 presents the corresponding optimal strategies $u_i^*(n)$ under different levels of the risk premium $\mu - r$. Specifically, we examine three scenarios: a low risk premium ($\mu - r = 0.03$), a moderate one $(\mu - r = 0.13)$, and a high one $(\mu - r = 0.25)$. The plots reveal that a larger risk premium, which implies a higher return on the risky asset relative to the risk‑free rate, leads to a higher optimal allocation to the risky asset (i.e., an increase in $u_i^*(n)$). This increase is clearly reflected in the vertical axis of Figure 3, where the level of $u_i^*(n)$ rises with the risk premium. The intuition is that investors demand greater compensation for bearing risk and thus increase their exposure when the risk return trade-off is more favorable. Conversely, a small risk premium offers only limited extra return for risk taking, which reduces the attractiveness of the risky asset and results in more conservative optimal strategies. These findings align with the predictions of portfolio theory and present the role of the risk premium in dynamic investment decisions.

\begin{figure}[htbp]
\centering
\begin{minipage}{\linewidth} 
    \centering
    \includegraphics[width=0.7\linewidth]{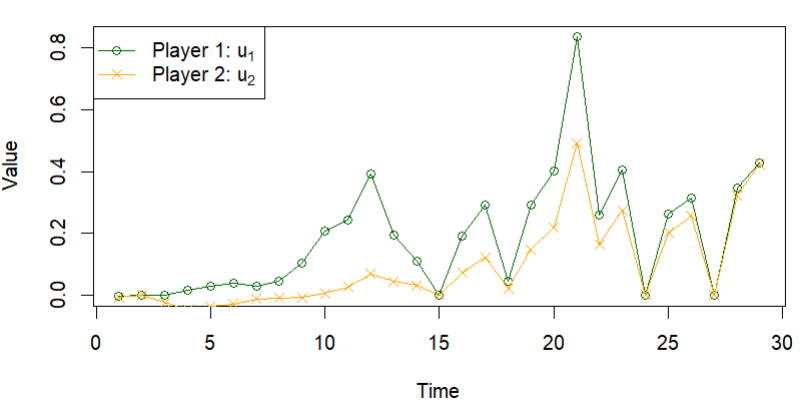} 
    \subfigure{ $(a)\quad \mu-r=0.03$} 
    \label{fig:group1_a}
\end{minipage}

\vspace{0.2cm} 

\begin{minipage}{\linewidth}
    \centering
    \includegraphics[width=0.7\linewidth]{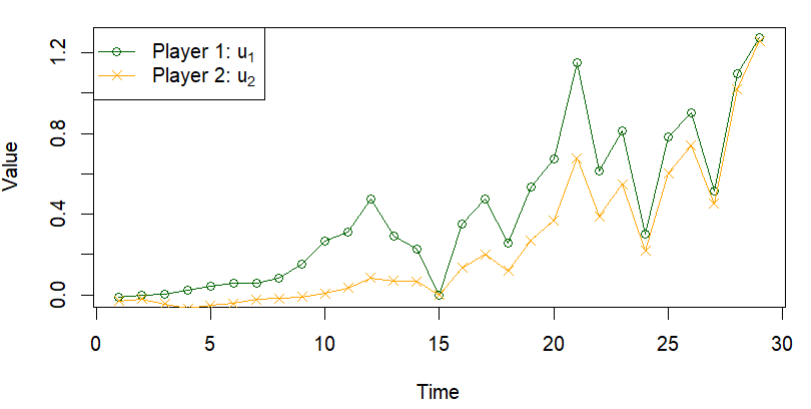}
    \subfigure{$(b)\quad \mu-r=0.13$} 
    \label{fig:group1_b}
\end{minipage}

\vspace{0.2cm}

\begin{minipage}{\linewidth}
    \centering
    \includegraphics[width=0.7\linewidth]{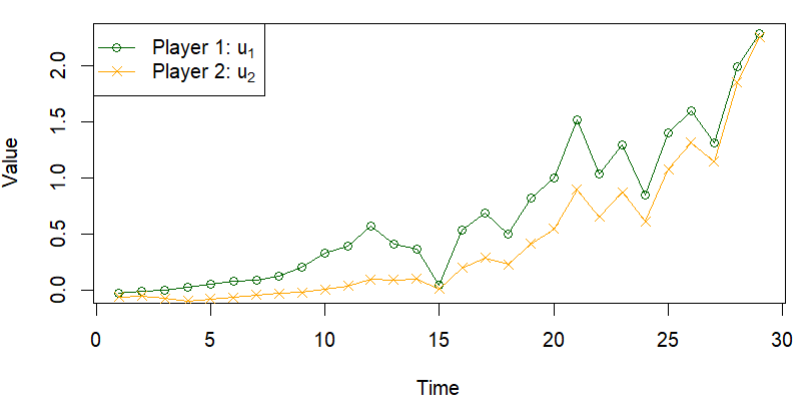}
    \subfigure{$(c)\quad \mu-r=0.25$} 
    \label{fig:group1_c}
\end{minipage}

\caption{Trajectories $u_1(n)$ and $u_2(n)$ under different $\mu - r$ dynamics}
\label{fig:group1}
\end{figure}
\section{Acknowledgments}

This work was supported by National Natural Science Foundation of China (Grant numbers 12471417) and National Key R$\&$D Program of China (Grant number 2023YFA1009200).

\end{document}